\documentclass{amsproc}
\usepackage{amssymb,amscd}
\usepackage[cp1251]{inputenc}
\usepackage[T2A]{fontenc}
\usepackage[matrix,arrow]{xy}

\newtheorem{theorem}{Theorem}
\newtheorem{definition}{Definition}
\newtheorem{lemma}{Lemma}
\newtheorem{proposition}{Proposition}
\newtheorem{remark}{Remark}

\newcommand{\bideg}{\mathrm{bideg}}

\newcommand{\supp}{\mathrm{supp}}

\newcommand{\Id}{\mathrm{Id}}

\begin{document}
\title[Mixed Hodge structure on complements of arrangements]{Mixed Hodge structure on complements of complex coordinate subspace arrangements}

\author{Yury~V. Eliyashev}
\address{Faculty of Mathematics, National Research University Higher School of Economics,  7 Vavilova Str., 117312 Moscow,
Russia} \address{Institute of Mathematics and Computer Science, Siberian Federal University, 79 Svobodny pr., 660041 Krasnoyarsk, Russia} \email{eliyashev@gmail.com}
\thanks{The work on Section 1-2
was supported by RScF grant, project 14-21-00053 dated 11.08.14.
 Section 3
was prepared within the framework of a subsidy granted to the HSE by the Government of the Russian Federation for the implementation of the Global Competitiveness Program. 
The work on Section 3 was supported by the grant of the Russian Federation Government to support scientific research under the supervision of a leading scientist at Siberian Federal University,  
№14.Y26.31.0006. The author is supported in part by a Dynasty Foundation grant.}
\begin{abstract}
We compute the mixed Hodge structure on the cohomology ring of complements
of complex coordinate subspace arrangements. 
The mixed Hodge structure can be described in terms of the special bigrading on the cohomology ring   of complements of complex coordinate subspace arrangements.
Originally this bigrading was introduced in the setting of toric topology by V.M. Buchstaber and T.E. Panov.
\end{abstract}

\maketitle

\noindent{\footnotesize\emph{Mathematics Subject Classification (2010).}
		32J25
	}\\
\noindent{\footnotesize\emph{Keywords.} Mixed Hodge structure, subspace arrangements.}	

\section*{Introduction}
A study of topology of coordinate subspace arrangements appears in
different areas of mathematics: in toric topology and combinatorial
topology \cite{BP,BP2}, in the theory of toric varieties, where
complements to  coordinate subspace arrangements play the role of
homogeneous coordinate spaces \cite{CX,CX2}, in the theory of
integral representations of holomorphic functions in several complex
variables, where coordinate subspace arrangements play the role of
singular sets of integral representations kernels \cite{AJ,Sh}.

The universal combinatorial method for the computation of cohomology
groups of complements to \emph{arbitrary} subspace arrangements was
developed in the book of Goresky and Macpherson \cite{GM} (see also
\cite{VS}), but this method often leads to cumbersome computations.
In the study of toric topology, in particular, in works of
Buchstaber and Panov \cite{BP, BP2}, the method for the computation
of the cohomology of complements to \emph{coordinate} subspace
arrangements was developed. This method is simpler than the
universal method and allows to get some additional topological
information.

The main purpose of this article is to compute the mixed Hodge structure
on the cohomology rings of complements to complex coordinate
subspace arrangements. We will show that
 this mixed Hodge structure is described by means of a special bigrading on the cohomology rings of complements to complex coordinate subspace arrangements, which was introduced in \cite{BP,BP2}.
 This bigrading was obtained originally from the combinatorial and topological ideas.

The mixed Hodge structure on the cohomology rings of complements to arbitrary \emph{hyperplane} arrangements in the complex and  $\ell$-adic settings was studied in \cite{Kim}, \cite{Leh}, \cite{Shap}. 
The same problem in the case of complements to arbitrary \emph{affine plane} arrangements  was studied in \cite{Dl2}.

The first section of this paper consists of different facts about
topology of complements to complex coordinate subspace arrangements.
In this section we follow \cite{BP}, \cite{BP2}.  Let
$Z_\mathcal{K}$ be a complex coordinate subspace arrangement in $\mathbb{C}^n.$ This arrangement is defined by the combinatorics of a simplicial complex $\mathcal{K}$ on the set $\{1,\dots,n\}.$
In \cite{BP}, \cite{BP2}, from the topological reasons, the
differential bigraded  algebra $R_\mathcal{K}$ was introduced ($R_\mathcal{K}$ is determined
by combinatorics of $\mathcal{K}$). The  cohomology ring
$H^*(\mathbb{C}^n \setminus Z_\mathcal{K})$ is isomorphic to the cohomology ring $H^*(R_\mathcal{K}).$ Denote by $H^{-p,2q}(R_\mathcal{K})$ the bigraded cohomology
of the algebra $R_\mathcal{K},$ then $$H^s(\mathbb{C}^n \setminus Z_\mathcal{K})\simeq
\bigoplus_{-p+2q=s} H^{-p,2q}(R_\mathcal{K}).$$ Thus, there is a bigrading on the
cohomology ring $H^*(\mathbb{C}^n \setminus Z_\mathcal{K})$.

In the second section we recall some facts and  concepts from
differential topology and algebraic geometry. We use them in
the last section.

In the third section the main theorem of this paper is proved. We
 show that the bigrading on the cohomology of $R_\mathcal{K}$ and,
consequently, the bigrading on the  cohomology
$H^*(\mathbb{C}^n\setminus Z_\mathcal{K})$ appear naturally from the mixed Hodge structure on cohomology of $\mathbb{C}^n\setminus Z_\mathcal{K}.$ In particular,
denote by $F^k H^s(\mathbb{C}^n\setminus Z_\mathcal{K},\mathbb{C})$ the $k$-th
term of the Hodge filtration on $H^s(\mathbb{C}^n\setminus
Z_\mathcal{K},\mathbb{C}),$ and denote by $W_k H^s(\mathbb{C}^n\setminus Z_\mathcal{K},\mathbb{C})$ the $k$-th
term of the weight filtration on $H^s(\mathbb{C}^n\setminus
Z_\mathcal{K},\mathbb{C}).$  Then following theorem holds.
\begin{theorem}$$F^{k} H^s(\mathbb{C}^n\setminus Z_{\mathcal{K}},\mathbb{C}) \cong \bigoplus_{\substack{q\geq k\\ -p +2 q =s}} H^{-p,2 q}(R_{\mathcal{K}})\otimes\mathbb{C},$$
$$W_{r} H^s(\mathbb{C}^n\setminus Z_{\mathcal{K}},\mathbb{C}) \cong \bigoplus_{\substack{2q\leq r\\ -p +2 q =s}} H^{-p,2 q}(R_{\mathcal{K}}) \otimes\mathbb{C}.$$
\end{theorem}

\section{General facts on topology of coordinate subspace arrangements}

This section consists of different facts about topology of complements of
coordinate subspaces arrangements. All statements of
this section are taken from \cite{BP}.

Let $\mathcal{K}$ be an arbitrary simplicial complex on the
set $[n]=\{1,\dots,n\},$ i.e., the vertices of $\mathcal{K}$ are elements of $[n].$ An element $j\in[n]$ is called a \emph{ghost vertex} 
if $j$ is not a vertex of $\mathcal{K}.$ If $\sigma=\{i_1,\dots,i_m\}$ is a subset of $[n],$ we write $\sigma\not\in \mathcal{K}$ if $\sigma$ does not define a simplex in $\mathcal{K}.$
Define the coordinate planes arrangement corresponding to $\mathcal{K}$ as
$$Z_{\mathcal{K}}:=\bigcup_{\sigma\not\in\mathcal{K}}L_\sigma,$$ where $$L_\sigma=\{z\in\mathbb{C}^n: z_{i}=0, i\in \sigma\}.$$
Every arrangement of complex coordinate subspaces in $\mathbb{C}^n$ can be defined in this way.

Let $D^2_\sigma\times S^1_\gamma$ be the chain:
$$D^2_\sigma\times S^1_\gamma=\{|z_i|\leq 1:i\in \sigma; |z_j|=1: j\in\gamma,z_k=1: k\not\in\gamma\cup\sigma\},$$
where $\sigma,\gamma\subseteq[n]$ and $\sigma\cap\gamma=\emptyset.$
Consider the differential form
 \begin{equation}\label{eq.dz}\frac{d z_I}{z_I}=\frac{d z_{i_1}}{z_{i_1}}\wedge\dots \wedge \frac{d z_{i_{k}}}{z_{i_{k}} },\end{equation}
 where $I\subseteq[n], |I|=k,$ $I=\{i_1,\dots,i_{k}\}$ and $i_1<\dots<i_{k}.$

 The orientation of the chain $D^2_\sigma\times S^1_\gamma$ is such that the restriction of the form
$$\frac{1}{(\sqrt{-1})^{|\gamma|}}\frac{d z_\gamma}{z_\gamma} \wedge \bigwedge_{j\in \sigma}(\sqrt{-1} d z_j\wedge d \overline{z}_j)$$
 to $D^2_\sigma\times S^1_\gamma$ is positive. Then the boundary of this chain equals
$$\partial D^2_\sigma\times S^1_\gamma=\Sigma_{i\in\sigma}(-1)^{(i,\gamma\cup i)+1}D^2_{\sigma\setminus i}\times S^1_{\gamma\cup i},$$
where $(i,\gamma\cup i)$ is the position of $i$ in the naturally ordered
set $\gamma\cup i.$

\begin{definition}\label{def.MA-comlex}The topological space
$$\mathcal{Z}_\mathcal{K}=\bigcup_{\sigma\in \mathcal{K}} D^2_\sigma\times S^1_{[n]\setminus \sigma}$$
is called the moment-angle complex.
\end{definition}

\begin{theorem}[\cite{BP}]\label{th.ret}
There exists a deformation retraction of $\mathbb{C}^n\setminus
Z_\mathcal{K}$ onto $\mathcal{Z}_\mathcal{K}.$
\end{theorem}

Let us introduce some formal algebraic notation. It will be used to describe the topology of moment-angle complexes.
\begin{definition}\label{def.SR-ring}
 A Stanley-Reisner ring of a simplicial complex  $\mathcal{K}$ on the vertex set $[n]$ is the ring
  $$\mathbb{Z}[\mathcal{K}]=\mathbb{Z}[v_1,\dots,v_n]/ \mathcal{I}_\mathcal{K},$$ where $\mathcal{I}_\mathcal{K}$ 
  is a homogeneous ideal generated by the monomials $v_\sigma=\prod_{i\in\sigma}v_i$ 
  such that $\sigma\not\in\mathcal{K}:$ $$\mathcal{I}_\mathcal{K}=(v_{i_1}\cdot{\dots}\cdot v_{i_m}: \{i_1,\dots,i_m\}\not\in\mathcal{K}).$$
\end{definition}

Consider the differential bigraded algebra
$(R(\mathcal{K}),\delta_R):$
$$R_\mathcal{K}:=\Lambda[u_1,\dots,u_n]\otimes
\mathbb{Z}[\mathcal{K}]/\mathcal{J},$$ where
$\Lambda[u_1,\dots,u_n]$ is the exterior algebra and $\mathcal{J}$ is
the ideal generated by monomials $v_i^2,u_i\otimes v_i,
i=1,\dots,n.$ The bidegrees of generators $v_i, u_i$  of this algebra
are equal to
$$\bideg\; v_i=(0,2), \bideg\; u_i=(-1,2).$$
The differential $\delta_R$ is defined on the generators as follows:
$$\delta_R u_i=v_i, \delta_R v_i =0.$$

Denote $u_I v_J:=u_{i_1}\dots u_{i_q}\otimes v_{j_1}\dots v_{j_p},$
where $I,J\subseteq[n],$ $I=\{i_1,\dots,i_q\},$
 $i_1<\dots<i_q,$ $J=\{j_1,\dots,j_p\},$ and $I\cap J=\emptyset$ (we suppose that $u_\emptyset v_\emptyset=1$).

Let $R^{-p,2q}_\mathcal{K}$ be the homogeneous component of bidegree $(-p,2q)$ of algebra $R_\mathcal{K}.$   The differential
$\delta_R$ is compatible with the bigrading, i.e.,
$\delta_R(R^{-p,2q}_\mathcal{K})\subseteq R^{-p+1,2q}_\mathcal{K}.$
Consider the complex
$$\dots\stackrel{\delta_R}{\longrightarrow}R^{-p-1,2q}_\mathcal{K}\stackrel{\delta_R}{\longrightarrow}R^{-p,2q}_\mathcal{K}\stackrel{\delta_R}{\longrightarrow}R^{-p+1,2q}_\mathcal{K}\stackrel{\delta_R}{\longrightarrow}\dots,$$ 
let $H^{-p,2q}(R_\mathcal{K})$ be a cohomology group of this
complex.  It is clear that the cohomology groups of $R_\mathcal{K}$ are
isomorphic to
$$H^s(R_\mathcal{K})=\bigoplus_{-p+2q=s}H^{-p,2q}(R_\mathcal{K}).$$

\begin{theorem}[\cite{BP}]\label{th.BP1}
The cohomology ring $H^*(\mathbb{C}^n\setminus Z_{\mathcal{K}})$ is
isomorphic to the ring $H^*(R_\mathcal{K}).$
\end{theorem}

\begin{remark}
The relation between Theorem \ref{th.BP1} and the results of Goresky
and Macpherson \cite{GM} on cohomology of subspace arrangements is
described in \cite[Ch.~8]{BP}.
\end{remark}

Let us describe an explicit construction of the isomorphism of
Theorem \ref{th.BP1}. First we construct a cell decomposition of
$\mathcal{Z}_\mathcal{K}.$ Consider the cell
$$E_{\sigma\gamma}=\{|z_i|<1:i\in \sigma; |z_j|=1,z_j\neq1: j\in\gamma;z_k=1: k\not\in\gamma\cup\sigma\},$$
where $\sigma,\gamma\subseteq[n]$ and $\sigma\cap\gamma=\emptyset.$
The closure of this cell equals
$\overline{E_{\sigma\gamma}}=D^2_\sigma\times S^1_\gamma.$ The
orientation of $E_{\sigma\gamma}$  is defined by the orientation of
$D^2_\sigma\times S^1_\gamma.$ We obtain the cell decomposition
$$\mathcal{Z}_\mathcal{K} = \bigcup_{\sigma\in \mathcal{K}, \gamma\subseteq [n]\setminus \sigma}E_{\sigma\gamma}.$$

Let $C_*(\mathcal{Z}_\mathcal{K})$ be the celluar chains group of
this cell decomposition; denote by
$C^*(\mathcal{Z}_\mathcal{K})$ the  celluar cochains group. Let
$E_{\sigma\gamma}'$ be the cocell dual to the cell
$E_{\sigma\gamma},$ i.e., $E_{\sigma\gamma}'$ is a linear function
on  $C_*(\mathcal{Z}_\mathcal{K})$ such that $\langle
E_{\sigma\gamma}',E_{\sigma'\gamma'}\rangle=\delta^{\sigma\gamma}_{\sigma'\gamma'}$
(the Kronecker delta).

\begin{proposition}[\cite{BP}]\label{th.BP2}
The linear map $\phi:R_\mathcal{K}\rightarrow
C^*(\mathcal{Z}_\mathcal{K}),$ $\phi(v_\sigma u_\gamma) =
E_{\sigma\gamma}'$ is an isomorphism of differential bigraded
modules. In particular, there exists an isomorphism
$H^*(R_\mathcal{K})\stackrel{\phi}{\simeq}
H^*(\mathcal{Z}_\mathcal{K}).$
\end{proposition}
From the structure of the cell decomposition of
$\mathcal{Z}_\mathcal{K}$ and Theorem \ref{th.ret} it follows that
every cycle $\Gamma\in H_s(\mathbb{C}^n\setminus Z_\mathcal{K})$ has
a representative of the form
\begin{equation}\label{eq.cycle_decom}
\Gamma=\sum_{-p+2q=s}\Gamma_{-p,2q},
\end{equation}
where $\Gamma_{-p,2q}$ is a cycle of the form
\begin{equation}\label{eq.cycle_type}\Gamma_{-p,2q}=\sum_{\substack{|\sigma|=q-p \\|\gamma|=p}}C_{\sigma \gamma}\cdot D^2_\sigma\times S^1_\gamma, C_{\sigma \gamma} \in \mathbb{Z}.\end{equation}
 Let $H_{-p,2q}(\mathbb{C}^n\setminus Z_{\mathcal{K}})$ be the group generated by all cycles of the form  (\ref{eq.cycle_type}).    Obviously, we have  $$H_{s}(\mathbb{C}^n\setminus Z_{\mathcal{K}})=\bigoplus_{-p+2q=s}H_{-p,2q}(\mathbb{C}^n\setminus Z_{\mathcal{K}}).$$

\begin{proposition}\label{pr.paribideg}
The pairing between the vector spaces
$H_{-p,2q}(\mathbb{C}^n\setminus Z_{\mathcal{K}},\mathbb{R})\subset H_{-p+2q}(\mathbb{C}^n\setminus Z_{\mathcal{K}},\mathbb{R})$ and
$\phi(H^{-p',2q'}(R_\mathcal{K}\otimes\mathbb{R}))\subset H^{-p'+2q'}(\mathbb{C}^n\setminus Z_{\mathcal{K}},\mathbb{R})$ is nondegenerate
if $p=p',$ $q=q'$ and is identically zero otherwise.
\end{proposition}
{\it Proof.}
Suppose that $p'\neq p$ and $q'\neq
q;$ then it follows from Proposition \ref{th.BP2} that
$\langle\Gamma_{-p,2q},\phi(\omega^{-p',2q'})\rangle=0$ for every
$\Gamma_{-p,2q}\in H_{-p,2q}(\mathbb{C}^n\setminus~Z_{\mathcal{K}})$ and every
$\omega^{-p',2q'}\in H^{-p',2q'}(R_\mathcal{K}).$ Hence the pairing between $\Gamma_{-p,2q}$ and
$\phi(\omega^{-p',2q'})$ can be nonzero only if $p'=p,q'=q.$ Since the pairing between $$H_s(\mathbb{C}^n\setminus Z_{\mathcal{K}},\mathbb{R})=\bigoplus_{-p+2q=s}H_{-p,2q}(\mathbb{C}^n\setminus Z_{\mathcal{K}},\mathbb{R})$$
and
$$H^s(\mathbb{C}^n\setminus Z_{\mathcal{K}},\mathbb{R})= \bigoplus_{-p+2q=s} \phi(H^{-p,2q}(R_\mathcal{K}\otimes\mathbb{R}))$$
is nondegenerate, we obtain the statement of the proposition.

\section{Mixed Hodge structures and resolvents of cycles}
Here we recall some facts about mixed Hodge structures. For references, see \cite{Dl}, \cite{KK}, \cite{PS}, \cite{Voi}.

\begin{definition}  Let $H$ be a finite-dimensional vector space over $\mathbb{Q}$. A pure Hodge structure of weight $s$ on $H$ is a decreasing filtration $F$ 
on  $H_\mathbb{C}=H \otimes \mathbb{C}$ such that $$F^p H_\mathbb{C} \cap\overline{ F^q H_\mathbb{C}}= 0$$
 whenever $p+q=s+1.$ The filtration $F$ is called the Hodge filtration.
\end{definition}

\begin{definition} Let $H$ be a finite-dimensional vector space over $\mathbb{Q}$. A mixed Hodge structure on $H$ consists, by definition, of the following:
\begin{enumerate}
  \item An increasing (weight) filtration $W$ on $H;$
  \item A decreasing (Hodge) filtration $F$ on  $H_\mathbb{С}=H \otimes \mathbb{C}$ satisfying the following condition: the filtration $F$ induces a pure Hodge structure of weight $s$ on $\mathrm{Gr}^W_s H= W_s H \otimes \mathbb{C} / W_{s-1} H \otimes \mathbb{C}.$
\end{enumerate}
\end{definition}


Cohomology groups of quasi-projective varieties admit a natural mixed Hodge structure. In particular,
 \begin{proposition} The weight filtration on the cohomology $H^s(X,\mathbb{Q})$ of a smooth variety $X$ has the form
 $$0=W_{s-1}\subset W_{s}\subset\dots \subset W_{2s}=H^s(X,\mathbb{Q}). $$
\end{proposition}

In this article we will use methods of differential topology, so we will consider the weight filtration on $H^s(X,\mathbb{C})$, not on $H^s(X,\mathbb{Q}).$ Now we recall the construction of the mixed Hodge structure on cohomology groups of a smooth complex algebraic variety.

Let $X$ be a smooth complex algebraic variety of dimension $n$. A
proper compactification of a variety $X$ is an open embedding
$j:X\hookrightarrow \overline{X}$ into a complete smooth algebraic
variety $\overline{X}$ such that $\overline{X}\setminus X= V
$ is a smooth normal crossing divisor. According
to Hironaka's theorem, a proper compactification always exists.

Let the divisor $V$ be defined by equations $z_1\cdot {\dots} \cdot z_k=0$ in a neighborhood $U\subset \overline{X},$ where $z_i$ are local coordinates in $U.$

\begin{definition} The sheaf
$$\Omega^m_{\overline{X}}(\log V)=\bigwedge^m (\Omega^1_{\overline{X}}(\log V))$$
is called the sheaf of holomorphic $m$-forms on $\overline{X}$ with
logarithmic poles along $V$, where $\Omega^1_{\overline{X}}$ is the
locally free $\mathcal{O}_{\overline{X}}$-module generated over $U$
by the differentials
$$\frac{d z_1}{z_1},\dots,\frac{d z_k}{z_k},d z_{k+1},\dots, d z_n.$$
\end{definition}

In other words, the sections of the sheaf
$\Omega^m_{\overline{X}}(\log V)$ in the neighborhood $U$ are
$m$-forms
$$\omega\wedge \frac{d z_I}{ z_I},$$
where $\omega$ is a holomorphic form on $U$, $I = \{i_1,\dots ,i_p\} \subset \{1,\dots ,k\}.$ Consider the following sheaves: $$\mathcal{E}^{p,q}_{\overline{X}}(\log V)=
\Omega^p_{\overline{X}}(\log V)
\otimes_{\mathcal{O}_{\overline{X}}}\mathcal{E}^{0,q}_{\overline{X}},$$  $$\mathcal{E}^{s}_{\overline{X}}(\log V)=\bigoplus_{p+q=s}\mathcal{E}^{p,q}_{\overline{X}}(\log V).$$
One can show that $$H^s(X,\mathbb{C})=H^s_d(\Gamma(\overline{X},\mathcal{E}^{s}_{\overline{X}}(\log V))),$$
where $\Gamma(\overline{X},\mathcal{E}^{s}_{\overline{X}}(\log V))$ is the space of global sections of  $\mathcal{E}^{s}_{\overline{X}}(\log V)$ and \linebreak 
$H^s_d(\cdot)$ is the cohomology with respect to the external derivative $d$.

Let us define the increasing weight filtration $W$ on the sheaf
$\Omega^m_{\overline{X}}(\log V)$  by setting
$$W_{k} \Omega^m_{\overline{X}}(\log V)= \Omega^{k}_{\overline{X}}(\log V)\wedge \Omega^{m-k}_{\overline{X}},$$
similarly on $\mathcal{E}^{p,q}_{\overline{X}}(\log V)$ and $\mathcal{E}^{s}_{\overline{X}}(\log V).$ Define the weight filtration on cohomology as follows
$$W_{k+s} H^s(X,\mathbb{C})=H^s_d(\Gamma(\overline{X},W_k \mathcal{E}^*_{\overline{X}}(\log V))).$$

Consider the decreasing Hodge filtration $$F^k
\mathcal{E}^s_{\overline{X}}(\log V)=\bigoplus_{p\geq k}\mathcal{E}^{p,s-p}_{\overline{X}}(\log V)$$ on
the sheaf $\mathcal{E}^s_{\overline{X}}(\log V).$ It induces the following Hodge filtration on the cohomology: $$F^k H^s(X,\mathbb{C})= H^s_d(\Gamma(\overline{X},F^k \mathcal{E}^*_{\overline{X}}(\log V))).$$
 \begin{proposition}\label{pr.logDR} These two filtration define the desired mixed Hodge structure on $H^s(X,\mathbb{C})$.\end{proposition}


Now we are going to develop some technical methods. We will use them to prove the main result of the paper. 
The main idea of what we are going to do is to define weight and Hodge filtrations on a kind of \v{C}ech-de Rham complex. Probably, the constructions below are well-known, but we don't know references.

Let $D_\alpha,\alpha\in \mathcal{A}$ be a finite set of divisors in $\overline{X}$ such that
$\bigcup_{\alpha\in \mathcal{A}}D_{\alpha}\cup V$ is a smooth normal
crossing divisor and $\bigcap_{\alpha\in \mathcal{A}}D_\alpha=\emptyset$.
This set of divisors defines an open cover $\mathcal{U}=\{\mathcal{U}_\alpha =\overline{X}\setminus D_\alpha\}_{\alpha\in\mathcal{A}}$ of $\overline{X}.$ 
Besides, it defines the open cover $\mathcal{U}_X=\{{\mathcal{U}_X}_\alpha =X\setminus D_\alpha\}_{\alpha\in\mathcal{A}}$ of $X.$

Consider the \v{C}ech complex $(C^{*}(\mathcal{U},\mathcal{E}^{r}(\log V)),\delta)$ of the cover $\mathcal{U}$ with coefficients in the sheaf $\mathcal{E}^{r}(\log V)$
$$\dots \stackrel{\delta}{\rightarrow} C^{t-1}(\mathcal{U},\mathcal{E}^{r}(\log V)) \stackrel{\delta}{\rightarrow} C^t(\mathcal{U},\mathcal{E}^{r}(\log V)) \stackrel{\delta}{\rightarrow} C^{t+1}(\mathcal{U},\mathcal{E}^{r}(\log V))\rightarrow\dots .$$
Let $(C^*(\mathcal{U},\mathcal{E}^{r}_{\log}),\delta)$ be a subcomplex of $(C^{*}(\mathcal{U},\mathcal{E}^{r}(\log V)),\delta)$ such that
for any $\omega\in C^t(\mathcal{U},\mathcal{E}^{r}_{\log})$ the element
$(\omega)_{\alpha_0,\dots,\alpha_t}$ is an element of
$\Gamma(\overline{X},\mathcal{E}^{r}_{\overline{X}}(\log V\cup
D_{\alpha_o}\cup\dots\cup D_{\alpha_t})).$
There exists a de Rham differential
$d:C^t(\mathcal{U},\mathcal{E}^{r}_{\log})\rightarrow
C^t(\mathcal{U},\mathcal{E}^{r+1}_{\log}),$
$$(d \omega)_{\alpha_0,\dots, \alpha_{t}}=d (\omega)_{\alpha_0,\dots, \alpha_{t}}.$$
The groups $C^t(\mathcal{U},\mathcal{E}^{r}_{\log})$ with the
differentials $d$ and $\delta$ form a double complex; its associated single complex is denoted by
$$K^s(\mathcal{U},\mathcal{E}_{\log})=\bigoplus_{r+t=s}C^t(\mathcal{U},\mathcal{E}^{r}_{\log}).$$
Consider the operator $D=(-1)^r\delta+d$ on $C^t(\mathcal{U},\mathcal{E}^{r}_{\log}).$ This operator defines a differential of the complex $K^s(\mathcal{U},\mathcal{E}_{\log})$. 
Hodge filtration $F$ and weight filtration $W$ are defined
naturally on the complexes $K^s(\mathcal{U},\mathcal{E}_{\log})$ and
$C^t(\mathcal{U},\mathcal{E}^{r}_{\log}).$

Consider a linear map $\varepsilon: \mathcal{E}^{s}_{\overline{X}}(\log V)\rightarrow C^0(\mathcal{U},\mathcal{E}^{s}_{\log})$ 
such that $\varepsilon(\omega)_\alpha=\omega|_{\mathcal{U}_\alpha},$ $\alpha\in \mathcal{A}.$ The map $\varepsilon: \mathcal{E}^{s}_{\overline{X}}(\log V)\rightarrow C^0(\mathcal{U},\mathcal{E}^{s}_{\log})$ 
induces a map from $\mathcal{E}^{s}_{\overline{X}}(\log V)$ to $K^s(\mathcal{U},\mathcal{E}_{\log}).$ We will denote the latter map by the same symbol $\varepsilon.$

Let $\rho_\alpha,\alpha\in \mathcal{A},$ be a partition of unity subordinated to the open cover $\mathcal{U}$, i.e., a set of real nonegative $C^\infty$-functions on $\overline{X}$ such that $\sum_{\alpha\in \mathcal{A}}\rho_\alpha\equiv1$ and $\supp(\rho_\alpha)\subset \mathcal{U}_\alpha.$
Using this partition of unity we define a homotopy operator $T: C^t(\mathcal{U},\mathcal{E}^{r}_{\log})\rightarrow C^{t-1}(\mathcal{U},\mathcal{E}^{r}_{\log})$:
$$(T\omega)_{i_0,\dots,i_{t-1}}=\sum_{\alpha\in \mathcal{A}}\rho_\alpha \omega_{\alpha,i_0,\dots,i_{t-1}};$$
here, $\omega$ is an element of $C^t(\mathcal{U},\mathcal{E}^{s}_{\log}).$ It is easy to check that
$$T\delta+ \delta T = \Id.$$

\begin{proposition}\label{pr.ChechDR}
The map $\varepsilon: \mathcal{E}^{*}_{\overline{X}}(\log V)\rightarrow K^*(\mathcal{U},\mathcal{E}_{\log})$ is a quasi-isomorphism of complexes. The induced isomorphism
 $$H^*(X,\mathbb{C}) \cong H^*(\mathcal{E}^{*}_{\overline{X}}(\log V),d)\stackrel{\varepsilon}{\cong} H^*(K^*(\mathcal{U},\mathcal{E}_{\log}),D)$$ respects Hodge and weight filtrations.
\end{proposition}
{\it Proof.}  Consider  the sequence
$$0\rightarrow \mathcal{E}^{s}_{\overline{X}}(\log V)\stackrel{\varepsilon}{\rightarrow}C^0(\mathcal{U},\mathcal{E}^{s}_{\log})\stackrel{\delta}{\rightarrow}C^1(\mathcal{U},\mathcal{E}^{s}_{\log})\stackrel{\delta}{\rightarrow}\dots$$
Since there exists a homotopy operator $T,$ this sequence is exact. Therefore, by standard arguments, $\varepsilon$ defines a quasi-isomorphism of complexes.
Moreover,  the homotopy operator is compatible with  weight and Hodge filtrations. In other words, we have $T(W_k C^t(\mathcal{U},\mathcal{E}^{s}_{\log})) \subset W_k C^{t-1}(\mathcal{U},\mathcal{E}^{s}_{\log})$ and $T(F^k C^t(\mathcal{U},\mathcal{E}^{s}_{\log})) \subset F^k C^{t-1}(\mathcal{U},\mathcal{E}^{s}_{\log}).$ 
Hence the sequences
$$0\rightarrow W_k\mathcal{E}^{s}_{\overline{X}}(\log V)\stackrel{\varepsilon}{\rightarrow}W_k C^0(\mathcal{U},\mathcal{E}^{s}_{\log})\stackrel{\delta}{\rightarrow} W_k C^1(\mathcal{U},\mathcal{E}^{s}_{\log})\stackrel{\delta}{\rightarrow}\dots,$$
$$0\rightarrow F^k\mathcal{E}^{s}_{\overline{X}}(\log V)\stackrel{\varepsilon}{\rightarrow}F^k C^0(\mathcal{U},\mathcal{E}^{s}_{\log})\stackrel{\delta}{\rightarrow} F^k C^1(\mathcal{U},\mathcal{E}^{s}_{\log})\stackrel{\delta}{\rightarrow}\dots$$
are exact and the induced filtrations on the cohomology groups of $K^*(\mathcal{U},\mathcal{E}_{\log})$ coincide with the filtrations on  $ H^*(\mathcal{E}^{*}_{\overline{X}}(\log V),d)$.

From now until the end of this section we will follow the ideas the paper
\cite{Gleason}. The main purpose of this part is to give definitions of a $\mathcal{U}_X$-chain and of a $\mathcal{U}_X$-resolvent of a cycle, and to describe their properties.
\begin{definition}
A $\mathcal{U}_X$-chain of degree $t$ and of dimension $s$ on the
variety ${X}$ is an alternating function $\Gamma$ from the set of
indexes $\mathcal{A}^{t+1}$ to the group of  $s$-dimensional singular chains of $X$
such that $\Gamma$ is nonzero at a finite number of
points from $\mathcal{A}^{t+1}$ and
$$\supp(\Gamma_{i_0,\dots,i_t})\subset
{\mathcal{U}_X}_{i_0}\cap\dots\cap{\mathcal{U}_X}_{i_t},$$  for every
$(i_0,\dots,i_{t})\in\mathcal{A}^{t+1},$ where
$\supp(\Gamma_{i_0,\dots,i_t})$ is the support of the chain
$\Gamma_{i_0,\dots,i_t}.$
\end{definition}

Let  $C_{t,s}({\mathcal{U}_X})$ be the additive group of
${\mathcal{U}_X}$-chains of degree $t$ and of dimension $s$ on the
variety ${X}.$ Define the maps $\delta':C_{t,s}(\mathcal{U}_X)\rightarrow
C_{t-1,s}(\mathcal{U}_X)$
$$(\delta'\Gamma)_{i_0,\dots,i_{t-1}}= \sum_{i\in\mathcal{A}}\Gamma_{i,i_0,\dots,i_{t-1}},$$
and $\partial:C_{t,s}(\mathcal{U}_X)\rightarrow
C_{t,s-1}(\mathcal{U}_X)$
$$(\partial\Gamma)_{i_0,\dots,i_{t}}=\partial(\Gamma)_{i_0,\dots,i_{t}},$$
i.e., the operator $\partial$ is a boundary operator on each chain
$\Gamma_{i_0,\dots,i_{t}}.$ 
Obviously, we have $\partial \partial =0$ and $\delta' \delta' =0.$ The groups $C_{t,s}(\mathcal{U}_X),t,s\geq
0,$ with the differentials $\delta',\partial$ form a double
complex. Let us define the map  $\varepsilon:C_{0,s}(\mathcal{U}_X)\rightarrow
C_{s}({X})$  as
$$\varepsilon(\Gamma)=\sum_{i\in\mathcal{A}}\Gamma_i,$$
where $C_{s}({X})$ is the group of the singular chains of dimension $s$ on $X$.

Now we will construct a pairing between elements of
$C_{t,s}(\mathcal{U}_X)$ and $C^t(\mathcal{U},\mathcal{E}^{s}_{\log}).$ Suppose that
$\Gamma \in C_{t,s}(\mathcal{U}_X)$ and $\omega\in
C^t(\mathcal{U},\mathcal{E}^{s}_{\log}).$ Then
$$\langle\omega,\Gamma\rangle= \frac{1}{(t+1)!}\sum_{(i_0,\dots,i_t)\in \mathcal{A}^{t+1}}\int_{\Gamma_{i_0,\dots,i_t}}\omega_{i_0,\dots,i_t}.$$

There exist the following relations for the pairing:
\begin{gather*}\label{eq.pairing}
\langle\omega^{t,s},\partial\Gamma_{t,s+1}\rangle=\langle d\omega^{t,s},\Gamma_{t,s+1}\rangle,\\
\langle\delta\omega^{t,s},\Gamma_{t+1,s}\rangle=\langle \omega^{t,s},\delta'\Gamma_{t+1,s}\rangle,\\
\int_{\varepsilon(\Gamma_{0,s})}\omega^s=\langle\varepsilon\omega^s,\Gamma_{0,s}\rangle,
\end{gather*}
where $\omega^{t,s}\in C^t(\mathcal{E}^s,\mathcal{U})$, $\omega^s\in
\mathcal{E}^s(X),$ and $\Gamma_{t,s} \in C_{t,s}(\mathcal{U}).$

\begin{definition}
Let $\Gamma$ be a singular cycle of dimension $s$ on ${X}.$ A
$\mathcal{U}_X$-resolvent of length $k$ of the cycle $\Gamma$ is a
collection of  $\mathcal{U}_X$-chains $\Gamma^i\in
C_{i,s-i}(\mathcal{U}_X),$ $0\leq~i\leq~k$ such that
$\Gamma=\varepsilon \Gamma^0$ and $\partial \Gamma^i=(-1)^{s-i}\delta'
\Gamma^{i+1}.$
\end{definition}

\begin{proposition}\label{prop.pairing}
Suppose that one has an $s$-dimensional cycle $\Gamma,$ a closed differential form
$\omega$ of degree $s$ on $X,$ a $\mathcal{U}_X$-resolvent
$\Gamma^0,\dots,\Gamma^k$ of the cycle $\Gamma,$ and a cocycle
$\widetilde{\omega}\in K^s(\mathcal{U},\mathcal{E}_{\log}).$ If
$\widetilde{\omega}=\sum_{i\leq k}{\widetilde{\omega}^{i,s-i}},$
$\widetilde{\omega}^{i,s-i}\in C^i(\mathcal{U},\mathcal{E}^{s-i}_{\log}),$
and the cocycle $\varepsilon \omega$ equals
$\widetilde{\omega}$ in
$H^s(K^*(\mathcal{U},\mathcal{E}_{\log}),D),$ then
$$\int_{\Gamma}\omega=\sum_{i\leq k}\langle\widetilde{\omega}^{i,s-i},\Gamma^i\rangle.$$
\end{proposition}
This proposition follows directly from the properties of the
pairing.
\begin{remark}
We considered $\mathcal{U}_X$-chains on the manifold ${X}$ with some special cover $\mathcal{U}_X$ and the pairing with elements of $C^t(\mathcal{U},\mathcal{E}^{s}_{\log}),$ because it is what we need in the sequel. 
However, one can consider the $\mathcal{U}_X$-chains for an arbitrary cover of a smooth manifold $X$ and the pairing of this $\mathcal{U}_X$-chains with the \v{C}ech cochains  $C^t(\mathcal{U}_X,\mathcal{E}^{s}).$ 
In this case the results above are also true.
\end{remark}

\section{Mixed Hodge structure on cohomology of complements to coordinate subspace arrangements}
In this section we compute mixed Hodge structure on the
cohomology ring
$H^*(\mathbb{C}^n\setminus~Z_{\mathcal{K}},\mathbb{C}).$ It follows
from Theorem \ref{th.BP1} and Proposition \ref{th.BP2} that there is
an isomorphism $H^*(\mathbb{C}^n\setminus
Z_{\mathcal{K}},\mathbb{C})\stackrel{\phi}{\cong}H^*(R_{\mathcal{K}}\otimes\mathbb{C}).$
Let us recall that we use the following notation: 
$H^{-p,2q}(R_{\mathcal{K}})$ is the bigraded cohomology group of the complex $R^{-p,2q}_{\mathcal{K}},$
$F^{k} H^s(\mathbb{C}^n\setminus Z_{\mathcal{K}},\mathbb{C})$ is the $k$-th term of the Hodge filtration on $H^s(\mathbb{C}^n\setminus Z_{\mathcal{K}},\mathbb{C}),$
$W_{r} H^s(\mathbb{C}^n\setminus Z_{\mathcal{K}},\mathbb{C})$ is the $r$-th term of the weight filtration on $H^s(\mathbb{C}^n\setminus Z_{\mathcal{K}},\mathbb{C}).$
\begin{theorem}\label{th.main}There exist the following isomorphisms: 
$$F^{k} H^s(\mathbb{C}^n\setminus Z_{\mathcal{K}},\mathbb{C}) \stackrel{\phi}{\cong} \bigoplus_{\substack{q\geq k\\ -p +2 q =s}} H^{-p,2 q}(R_{\mathcal{K}}\otimes\mathbb{C}),$$
$$W_{r} H^s(\mathbb{C}^n\setminus Z_{\mathcal{K}},\mathbb{C}) \stackrel{\phi}{\cong} \bigoplus_{\substack{2q\leq r\\ -p +2 q =s}} H^{-p,2 q}(R_{\mathcal{K}}\otimes\mathbb{C}).$$
\end{theorem}
{\it Proof.} Let us construct a proper compactification of $\mathbb{C}^n\setminus
Z_{\mathcal{K}}.$ Consider the standard embedding
$\mathbb{C}^n\setminus Z_{\mathcal{K}}\hookrightarrow \mathbb{C}^n
\hookrightarrow \mathbb{CP}^n.$ Then
$$\mathbb{C}^n\setminus Z_{\mathcal{K}} = \mathbb{CP}^n \setminus (\mathbb{CP}^{n-1}_{\infty}\cup\bigcup_{\sigma\not\in \mathcal{K}} \overline{L_{\sigma}}),$$
where $\overline{L_{\sigma}}$ is the closure of the complex plane $$L_\sigma=\{z\in\mathbb{C}^n: z_{i}=0, i\in \sigma\},$$
and $\mathbb{CP}^{n-1}_{\infty}$ is the hyperplane at infinity. Making a sequence of blowups along irreducible components of $\mathbb{CP}^{n-1}_{\infty}\cup\bigcup_{\sigma\not\in \mathcal{K}} \overline{L_{\sigma}}$
we get a proper compactification of $\mathbb{C}^n\setminus Z_{\mathcal{K}}.$ Denote this compactification by $\overline{X}.$ 
Let $\widetilde{L}_i$ be the proper preimage of the closure of $L_i=\{z\in\mathbb{C}^n: z_i=0\}$ in $\overline{X}.$ 
Put $$\mathcal{U}_\sigma=\overline{X}\setminus \bigcup_{i\not \in \sigma}\widetilde{L}_i.$$ Then
$\mathcal{U}_{\mathcal{K}}=\{\mathcal{U}_\sigma\}_{\sigma \in \mathcal{K}}$ is the open cover of $\overline{X}.$ 
The restriction of $\mathcal{U}_{\mathcal{K}}$ to $\mathbb{C}^n\setminus Z_{\mathcal{K}}$ is the open cover of $\mathbb{C}^n\setminus Z_{\mathcal{K}}.$ 
If $\widetilde{\mathcal{U}}_{\sigma}=\mathcal{U}_{\sigma} \cap (\mathbb{C}^n\setminus Z_{\mathcal{K}}),$ then $\widetilde{\mathcal{U}}_{\sigma}$ is isomorphic to $\mathbb{C}^{|\sigma|}\times (\mathbb{C}^*)^{n-|\sigma|}.$

Consider the logarithmic \v{C}ech-de Rham double complex
$(C^t(\mathcal{U}_{\mathcal{K}},\mathcal{E}^{r}_{\log}),d,\delta)$ and the associated
single complex $(K^s(\mathcal{U}_{\mathcal{K}},\mathcal{E}_{\log}),D).$ It follows from Proposition \ref{pr.ChechDR} that cohomology of $H^s(K^*(\mathcal{U}_{\mathcal{K}},\mathcal{E}_{\log}),D)$ is
isomorphic to the de Rham cohomology group
$H^s(\mathbb{C}^n\setminus Z_\mathcal{K},\mathbb{C}).$ Weight filtrations on both cohomology groups
coincide, and the same is true for Hodge filtration. Define a subcomplex
$(M^{r,t}(\mathcal{U}_\mathcal{K}),d,\delta)$ of
$(C^t(\mathcal{U},\mathcal{E}^{r}_{\log}),d,\delta);$ the elements of this
complex are $$(\omega)_{\alpha_0 \dots \alpha_t} =
\sum_{\substack{|I|=r\\I\subset [n]\setminus (\alpha_0\cap \dots\cap
\alpha_t)}} C^I_{\alpha_0 \dots \alpha_t} \frac{d z_I}{z_I},
$$ where $C^I_{\alpha_0 \dots \alpha_t} \in \mathbb{C}$ and $\alpha_0, \dots,
\alpha_t \in \mathcal{K}.$ It is easy
to see that $\frac{d z_I}{z_I}$ are indeed logarithmic forms in the
proper compactification. Denote by
$j:M^{r,t}(\mathcal{U}_\mathcal{K})\hookrightarrow
C^t(\mathcal{U},\mathcal{E}^{r}_{\log})$ the inclusion of the subcomplex.
The associated single complex is denoted
$M^{s}(\mathcal{U}_\mathcal{K}).$ Observe that the differential $d$ acts
trivially on the complex $M^{r,t}(\mathcal{U}_\mathcal{K}),$ whence
the cohomology $H^s(M^{\bullet}(\mathcal{U}_\mathcal{K}),D)$ is
isomorphic to
$\bigoplus_{r+t=s}H^t(M^{r,\bullet}(\mathcal{U}_\mathcal{K}),\delta),$
where $H^t(M^{r,\bullet}(\mathcal{U}_\mathcal{K}),\delta)$ is the
cohomology of the complex
$$\dots\rightarrow M^{r,t-1}(\mathcal{U}_\mathcal{K}) \rightarrow M^{r,t}(\mathcal{U}_\mathcal{K})\rightarrow M^{r,t+1}(\mathcal{U}_\mathcal{K})\rightarrow\dots.$$

\begin{lemma}
The inclusion map $j:M^{s}(\mathcal{U}_\mathcal{K})\hookrightarrow
K^{s}(\mathcal{U}_\mathcal{K},\mathcal{E}_{\log})$ is a quasi-isomorphism.
\end{lemma}
{\it Proof of the Lemma.} We will show that for any element
$\varphi\in K^{s}(\mathcal{U}_\mathcal{K},\mathcal{E}_{\log})$ such that $D \varphi=0,$
there exist $\psi\in M^{s}(\mathcal{U}_\mathcal{K})$ and $\omega\in
K^{s-1}(\mathcal{U}_\mathcal{K},\mathcal{E}_{\log})$ such that $D \psi=0,$ $\varphi =
\psi + D \omega.$ This will prove the lemma.
First, observe that any cocycle in $K^{s}(\mathcal{U}_\mathcal{K},\mathcal{E}_{\log})$ is cohomologous to a cocycle from $C^0(\mathcal{U}_\mathcal{K},\mathcal{E}^{s}_{\log}).$ 
Indeed, by Proposition \ref{pr.ChechDR}
the map $\varepsilon: \mathcal{E}^{s}_{\overline{X}}(\log V)\rightarrow C^0(\mathcal{U}_\mathcal{K},\mathcal{E}^{s}_{\log})\subset K^s(\mathcal{U}_\mathcal{K},\mathcal{E}_{\log})$ is a quasi-isomorphism of complexes.
Let $\varphi$ be a cocycle from $C^0(\mathcal{U}_\mathcal{K},\mathcal{E}^{s}_{\log}).$
Let us prove that for any $k\geq -1$ there exist $\psi^0,\dots,\psi^k,\psi^i\in M^{s-i,i}(\mathcal{U}_\mathcal{K}), D \psi^i=0,$ $\omega^0,\dots,\omega^k,\omega^i\in C^i(\mathcal{U}_\mathcal{K},\mathcal{E}^{s-1-i}_{\log}),$ and $\varphi^k\in C^{k+1}(\mathcal{U}_\mathcal{K},\mathcal{E}^{s-k-1}_{\log}), D \varphi^k =0$ such that
$$\varphi= \sum_{i=0}^k \psi^i + D \sum_{i=0}^k \omega^i + \varphi^k.$$

We will construct inductively such a cocycle decomposition.
The base of induction: $k=-1,$ this case is trivial, $\varphi= \varphi^{-1}$.

Suppose that the decomposition $$\varphi= \sum_{i=0}^k \psi^i + D \sum_{i=0}^k \omega^i + \varphi^k$$
 is already constructed for a given $k$.

 Since $D \varphi^k=0,$ the form $\varphi^k_{\sigma_0,\dots,\sigma_{k+1}},\sigma_0,\dots,\sigma_{k+1}\in \mathcal{K}$ is a closed $(s-k-1)$-form on $\widetilde{\mathcal{U}}_{\sigma_0\cap\dots\cap \sigma_{k+1}}$. 
 There exists a unique decomposition $$\varphi^k_{\sigma_0,\dots,\sigma_{k+1}}=\sum_{\substack{|I|=s-k-1\\I\subset [n]\setminus (\sigma_0\cap\dots\cap \sigma_{k+1})}}C^I_{\sigma_0,\dots,\sigma_{k+1}} \frac{d z_I}{z_I}+ d \omega_{\sigma_0,\dots,\sigma_{k+1}}.$$ 
 Indeed, differential forms $\frac{d z_I}{z_I},|I|=s-k-1,I\subset [n]\setminus (\sigma_0\cap\dots\cap \sigma_{k+1})$ are a basis of $H^{s-k-1}(\widetilde{\mathcal{U}}_{\sigma_0\cap\dots\cap \sigma_{k+1}}).$ 
 Putting $\omega^{k+1}_{\sigma_0,\dots,\sigma_{k+1}}=-\omega_{\sigma_0,\dots,\sigma_{k+1}},$ $\varphi^{k+1}=(-1)^{s-k-1} \delta \omega^{k+1},$ 
 and $$\psi^{k+1}_{\sigma_0,\dots,\sigma_{k+1}}=\sum_{\substack{|I|=s-k-1\\I\subset [n]\setminus (\sigma_0\cap\dots\cap \sigma_{k+1})}}C^I_{\sigma_0,\dots,\sigma_{k+1}} \frac{d z_I}{z_I},$$ 
 we see that
 $$\varphi= \sum_{i=0}^{k+1} \psi^i + D \sum_{i=0}^{k+1} \omega^i + \varphi^{k+1}.$$
 Let us check that $D \psi^{k+1}=0, D \varphi^{k+1}=0.$ From definition we have
 $$D \psi^{k+1}+ D \varphi^{k+1}=0, d \psi^{k+1}=0, \delta \varphi^{k+1} =0,$$
 therefore,  $(-1)^{s-k-1}\delta \psi^{k+1} = - d \psi^{k+1}$ and
 $(\delta \psi^{k+1})_{\sigma_0,\dots,\sigma_{k+2}}$ is a liner combination of differential forms $\frac{d z_I}{z_I}$ on $\widetilde{\mathcal{U}}_{\sigma_0\cap\dots\cap \sigma_{k+2}}\cong \mathbb{C}^{|\sigma_0\cap\dots\cap \sigma_{k+2}|}\times (\mathbb{C}^*)^{n-|\sigma_0\cap\dots\cap \sigma_{k+2}|}.$
 This form is exact if and only if it is identically zero. Thus, $(-1)^{s-k-1}\delta \psi^{k+1} = -d \psi^{k+1}=0.$ The induction step is proved.

 Taking $k=s$ we have
$$\varphi= \sum_{i=0}^s \psi^i + D \sum_{i=0}^s \omega^i + 0,$$
for $\psi = \sum_{i=0}^s \psi^i$ and $\omega = \sum_{i=0}^s \omega^i,$
we have
$$\varphi= \psi + D \omega,$$
and proves the lemma. \hfill$\square$

\begin{lemma}\label{lemma.filtr} Mixed Hodge structure on $H^s(K^*(\mathcal{U}_{\mathcal{K}},\mathcal{E}_{\log}),D) \cong H^s(M^{*}(\mathcal{U}_\mathcal{K}),D)$ has the form 
$$F^k H^s(M^{*}(\mathcal{U}_\mathcal{K}),D)=\bigoplus_{\substack{r\geq k\\r+t=s}}H^t(M^{r,\bullet}(\mathcal{U}_\mathcal{K}),\delta),$$
$$W_k H^s(M^{*}(\mathcal{U}_\mathcal{K}),D)=\bigoplus_{\substack{2r\leq k\\r+t=s}}H^t(M^{r,\bullet}(\mathcal{U}_\mathcal{K}),\delta).$$
\end{lemma}
{\it Proof of the Lemma.} The complex $(M^{r,t}(\mathcal{U}_\mathcal{K}),\delta)$ is naturally isomorphic to the complex $(C^{t}(\mathcal{U}_\mathcal{K},H^r(\bullet)),\delta).$ 
An element $\psi\in C^{t}(\mathcal{U}_\mathcal{K},H^r(\bullet))$ is a cochain with coefficients in cohomology, i.e., $\psi_{\sigma_0,\dots,\sigma_{t}}\in H^r(\widetilde{\mathcal{U}}_{\sigma_0}\cap\dots\cap \widetilde{\mathcal{U}}_{\sigma_{t}},\mathbb{C}).$ 
There exists a natural mixed Hodge structure on $H^r(\widetilde{\mathcal{U}}_{\sigma_0}\cap\dots\cap \widetilde{\mathcal{U}}_{\sigma_{t}},\mathbb{C}):$
$$H^*(\widetilde{\mathcal{U}}_{\sigma_0}\cap\dots\cap \widetilde{\mathcal{U}}_{\sigma_{t}},\mathbb{C})\simeq{\bigwedge}^*[\frac{d z_i}{z_i}:i\not\in {\sigma_0}\cap\dots\cap {\sigma_{t}}],$$
$$F^k H^*(\widetilde{\mathcal{U}}_{\sigma_0}\cap\dots\cap \widetilde{\mathcal{U}}_{\sigma_{t}},\mathbb{C})\simeq \bigoplus_{r\geq k} {\bigwedge}^r[\frac{d z_i}{z_i}:i\not\in {\sigma_0}\cap\dots\cap {\sigma_{t}}],$$
$$W_k H^*(\widetilde{\mathcal{U}}_{\sigma_0}\cap\dots\cap \widetilde{\mathcal{U}}_{\sigma_{t}},\mathbb{C})\simeq \bigoplus_{2r\leq k} {\bigwedge}^r[\frac{d z_i}{z_i}:i\not\in {\sigma_0}\cap\dots\cap {\sigma_{t}}].$$
By functoriality we obtain a mixed Hodge structure on $H^s(M^{*}(\mathcal{U}_\mathcal{K}),D).$ This gives us the statement of the lemma.\hfill$\square$

\begin{lemma}\label{lemma.resolv} Let $$\Gamma_{-p,2q}=\sum_{\substack{|\sigma|=q-p\\|\gamma|=p}}C^{\sigma \gamma}\cdot D^2_\sigma\times S^1_\gamma$$ be a cycle in $\mathbb{C}^n\setminus Z_{\mathcal{K}}.$ 
Then there exists a $\mathcal{U}_\mathcal{K}$-resolvent $\Gamma^{0}_{-p,2q},\dots,\Gamma^{q-p}_{-p,2q}$ of length $q-p,$ where $\Gamma^{k}_{-p,2q}$ is a $\mathcal{U}_\mathcal{K}$-chain of dimension $2q-p-k$ and of degree $k$ of the form
$$(\Gamma^{k}_{-p,2q})_{\alpha_0,\dots,\alpha_k}=\sum_{\substack{|\sigma|=q-p-k\\|\gamma|=p+k}}C^{\sigma \gamma}_{ \alpha_0\dots\alpha_k}\cdot D^2_\sigma\times S^1_\gamma.$$
\end{lemma}
{\it Proof of the Lemma.} We will use induction on the length $k$ of the
resolvent. We are going to construct a resolvent of the special form
$$(\Gamma^{k}_{-p,2q})_{\sigma_k,\sigma_{k-1},\dots,\sigma_0}=\sum_{|\gamma|=p+k}C^{\sigma_k \gamma}_{\sigma_k\dots \sigma_0}\cdot D^2_{\sigma_k}\times S^1_\gamma$$
for $|\sigma_j|=q-p-j,$  $\sigma_{j}\subset \sigma_t, j>t,$
$t,j=0,\dots,k,$ and
$(\Gamma^{k}_{-p,2q})_{\alpha_0,\dots,\alpha_k}=0$ for any other
indexes $\alpha_0,\dots,\alpha_k.$

The base of induction: define
 $(\Gamma^0_{-p,2q})_{\sigma_0}=\sum_{\substack{|\gamma|=p}}C^{\sigma_0 \gamma}\cdot D^2_{\sigma_0}\times S^1_\gamma$ with $|\sigma_0|=q-p,$
and $(\Gamma^0_{-p,2q})_{\alpha}=0$ for any other indexes $\alpha.$ We
get
$$\Gamma_{-p,2q}=\sum_{\substack{|\sigma_0|=q-p\\|\gamma|=p}}C^{\sigma_0 \gamma}\cdot D^2_{\sigma_0}\times S^1_\gamma=\sum_{\sigma\in\mathcal{K}}(\Gamma^{0}_{-p,2q})_{\sigma}=\varepsilon' \Gamma^{0}_{-p,2q},$$
therefore $\Gamma^{0}_{-p,2q}$ is a resolvent of length $0$.

Suppose that a resolvent $\Gamma^{0}_{-p,2q},\dots,\Gamma^{k}_{-p,2q}$
of length $k$ is already constructed. Recall that $(i,\gamma)$  is
the position of $i$ in the naturally ordered set $\gamma\cup i.$
Put $$(\Gamma^{k+1}_{-p,2q})_{\sigma_k\setminus i,\sigma_k
\dots\sigma_0} =(-1)^{2q-p-k}\sum_{\substack{|\gamma|=p+k}}
(-1)^{(i,\gamma)}C^{\sigma_k \gamma}_{ \sigma_k \dots\sigma_0}
D^2_{\sigma_k\setminus i}\times S^1_{\gamma\cup i}$$ for $i\in
\sigma_k,$ $|\sigma_j|=q-p-j,$ $\sigma_{j+1}\subset \sigma_{j},$ and put
$$(\Gamma^{k+1}_{-p,2q})_{\alpha_0,\dots,\alpha_{k+1}}=0$$ for any
other indexes $\alpha_0,\dots,\alpha_{k+1}.$ Let us show that
$\Gamma^{0}_{-p,2q},\dots,\Gamma^{k+1}_{-p,2q}$ is a resolvent of length
$k+1.$ We have
\begin{multline*}(-1)^{2 q - p -k}(\delta' \Gamma^{k+1}_{-p,2q})_{\sigma_k \dots\sigma_0}=\sum_{i\in \sigma_k}(\Gamma^{k+1}_{-p,2q})_{\sigma_k\setminus i,\sigma_k \dots\sigma_0}=\\=\sum_{\substack{i\in \sigma_k\\|\gamma|=p+k}} (-1)^{(i,\gamma)}C^{\sigma_k \gamma}_{ \sigma_k \dots\sigma_0} D^2_{\sigma_k\setminus i}\times S^1_{\gamma\cup i}
=\sum_{|\gamma|=p+k}C^{\sigma_k \gamma}_{ \sigma_k \dots\sigma_0}
\partial D^2_{\sigma_k}\times S^1_{\gamma}=\\=(\partial
\Gamma^{k}_{-p,2q})_{\sigma_k \dots\sigma_0}.\end{multline*}
 For any indexes $\alpha_0,\dots,\alpha_{k}$ different from $\sigma_k,\dots,\sigma_{m+1},\sigma_{m-1},\dots,\sigma_0,$ $0\leq m \leq k,$  directly from definition of the chain
 $\Gamma^k_{-p,2q},\Gamma^{k+1}_{-p,2q},$ we get $$(\partial\Gamma^{k}_{-p,2q})_{\alpha_0,\dots,\alpha_{k}}=(-1)^{2q-p-k}(\delta'\Gamma^{k+1}_{-p,2q})_{\alpha_0,\dots,\alpha_{k}}=0.$$
Consider the last case $\sigma_k\setminus
i,\sigma_k,\dots,\sigma_{m+1},\sigma_{m-1},\dots,\sigma_0$ for
$0\leq m \leq k.$ Since by the induction hypothesis
$\Gamma^{0}_{-p,2q}\dots\Gamma^{k}_{-p,2q}$  is a resolvent, we get
$$(-1)^{2q-p-k+1}\delta'\Gamma^{k}_{-p,2q}=\partial
\Gamma^{k-1}_{-p,2q},$$ whence $\delta' \partial
\Gamma^{k}_{-p,2q}=0$ and
\begin{multline*}(\delta' \partial\Gamma^{k}_{-p,2q})_{\sigma_k\dots\sigma_{m+1}\sigma_{m-1}\dots\sigma_0}=\\=
\sum_{\substack{\sigma_{m+1}\subset\sigma_m\subset\sigma_{m-1}\\|\sigma_m|=q-p-m}}\sum_{|\gamma|=p+k}\sum_{i\in \sigma_k} (-1)^{(i,\gamma)} C^{\sigma_k \gamma}_{ \sigma_k\dots\sigma_{0}}\cdot D^2_{\sigma_k\setminus i}\times S^1_{\gamma\cup i}=0,\end{multline*}
 Therefore, for a fixed $i\in
\sigma_k$ we get
$$\sum_{\substack{\sigma_{m+1}\subset\sigma_m\subset\sigma_{m-1}\\|\sigma_m|=q-p-m}}\sum_{|\gamma|=p+k} (-1)^{(i,\gamma)} C^{\sigma_k \gamma}_{ \sigma_k\dots\sigma_{0}}\cdot D^2_{\sigma_k\setminus i}\times S^1_{\gamma\cup i}=0.$$
On the other hand,
\begin{multline*}(\delta'\Gamma^{k+1}_{-p,2q})_{\sigma_k\setminus
i,\sigma_k
\dots\sigma_{m+1}\sigma_{m-1}\dots\sigma_0}=\\=\sum_{\substack{\sigma_{m+1}\subset\sigma_m\subset\sigma_{m-1}\\|\sigma_m|=q=p-m}}\sum_{|\gamma|=p+k}(-1)^{(i,\gamma)}
C^{\sigma_k \gamma}_{ \sigma_k\dots\sigma_{0}}\cdot
D^2_{\sigma_k\setminus i}\times S^1_{\gamma\cup i},\end{multline*}
whence $(\delta'\Gamma^{k+1}_{-p,2q})_{\sigma_k\setminus
i,\sigma_k \dots\sigma_{m+1}\sigma_{m-1}\dots\sigma_0}=0.$ We have shown that 
 $$\partial
\Gamma^{k}_{-p,2q}=(-1)^{2q-p-k}\delta' \Gamma^{k+1}_{-p,2q}.$$
\hfill
$\Box$

Recall that by Proposition \ref{pr.ChechDR}
 the map $\varepsilon: H^*(\mathbb{C}^n\setminus~Z_{\mathcal{K}},\mathbb{C})\rightarrow H^*( K^*(\mathcal{U}_\mathcal{K},\mathcal{E}_{\log}))$ is an isomorphism.
\begin{lemma}\label{lemma.PR} Let $$\Gamma_{-p,2q}=\sum_{\substack{|\sigma|=q-p\\|\gamma|=p}}C^{\sigma \gamma}\cdot D^2_\sigma\times S^1_\gamma$$ 
be a cycle in $\mathbb{C}^n\setminus Z_{\mathcal{K}},$ and let $\psi\in M^{r,t}(\mathcal{U}_\mathcal{K})\subseteq K^{r+t}(\mathcal{U}_\mathcal{K},\mathcal{E}_{\log})$ be a cocycle. 
Then $$\int_{\Gamma_{-p,2q}} \varepsilon^{-1}\psi = 0$$
for any $p,q,r,t$ such that $r\neq q$ or $t \neq q-p.$
\end{lemma}
{\it Proof of the Lemma.}
We may assume that $2q-p=r+t=s$, otherwise $\dim\Gamma_{-p,2q}\neq \deg \varepsilon^{-1}\psi$ and the integral is automatically zero.
By Lemma \ref{lemma.resolv} we have the resolvent  $\Gamma^{0}_{-p,2q},\dots,\Gamma^{q-p}_{-p,2q}$ of $\Gamma_{-p,2q}.$

The first case is $r>q.$  By Proposition \ref{prop.pairing} we have
$$\int_{\Gamma_{-p,2q}} \varepsilon^{-1}\psi = \langle \psi,\Gamma^{t}_{-p,2q}\rangle,$$
where
$$(\Gamma^{t}_{-p,2q})_{\alpha_0,\dots,\alpha_t}=\sum_{\substack{|\sigma|=q-p-t\\|\gamma|=p+t}}C^{\sigma \gamma}_{ \alpha_0\dots\alpha_t}\cdot D^2_\sigma\times S^1_\gamma.$$
Since $q-p-t>0,$
the pairing $\langle \psi,\Gamma^{t}_{-p,2q}\rangle$ is a sum of integrals of the differential forms $\frac{d z_I}{ z_I}$ over chains $D^2_\sigma\times S^1_\gamma, |\sigma|\neq0.$ 
A direct computation shows that all this integrals are equal to zero. So,
$$\int_{\Gamma_{-p,2q}} \varepsilon^{-1}\psi=0.$$

The second case is $r<q.$ There exist cochains $\omega^{r,s-r-1},\dots,\omega^{q-1,s-q},$ $\omega^{i,j}\in C^j(\mathcal{U}_{\mathcal{K}},\mathcal{E}^{i}_{\log}),$ 
such that the cocycle $$\varphi=\psi + D\sum_{i=r}^{q-1} \omega^{i,s-i-1}$$
is an element of $C^{s-q}(\mathcal{U}_{\mathcal{K}},\mathcal{E}^{q}_{\log}).$
By Proposition \ref{prop.pairing} we have
$$\int_{\Gamma_{-p,2q}} \varepsilon^{-1}\psi = \langle \varphi,\Gamma^{q-p}_{-p,2q}\rangle,$$
where
$$(\Gamma^{q-p}_{-p,2q})_{\alpha_0,\dots,\alpha_t}=\sum_{|\gamma|=q}C^{\gamma}_{ \alpha_0\dots\alpha_{q-p}}\cdot S^1_\gamma.$$
Observe that the differential forms $\varphi_{\sigma_0,\dots,\sigma_{s-q}}=d(\omega^{q-1,s-q}_{\sigma_0,\dots,\sigma_{s-q}})$ are exact.
Thus, the pairing $\langle \varphi,\Gamma^{q-p}_{-p,2q}\rangle$ is a sum of integrals of exact differential forms over cycles $S^1_\gamma,$ and all such integrals are equal to zero. So, we have
$$\int_{\Gamma_{-p,2q}} \varepsilon^{-1}\psi=0.$$
The lemma is proved. \hfill$\square$

By Proposition \ref{pr.paribideg} the pairing between the vector spaces
$H_{-p,2q}(\mathbb{C}^n\setminus Z_{\mathcal{K}},\mathbb{C})\subset H_{-p+2q}(\mathbb{C}^n\setminus Z_{\mathcal{K}},\mathbb{C})$ and
$\phi(H^{-p',2q'}(R_\mathcal{K}\otimes\mathbb{C}))\subset H^{-p'+2q'}(\mathbb{C}^n\setminus Z_{\mathcal{K}},\mathbb{C})$ is nondegenerate
if $p=p',$ $q=q',$ and is identically zero otherwise.
On the other hand,  $$H^*(\mathbb{C}^n\setminus~Z_{\mathcal{K}},\mathbb{C})\stackrel{\varepsilon}{\simeq} H^*( K^*(\mathcal{U}_\mathcal{K},\mathcal{E}_{\log}))\simeq H^s(M^{\bullet}(\mathcal{U}_\mathcal{K}),D)\simeq\bigoplus_{r+t=s}H^t(M^{r,\bullet}(\mathcal{U}_\mathcal{K}),\delta).$$ 
By Lemma \ref{lemma.PR} we see that the pairing between the spaces $H_{-p,2q}(\mathbb{C}^n\setminus Z_{\mathcal{K}},\mathbb{C})\subset H_{-p+2q}(\mathbb{C}^n\setminus Z_{\mathcal{K}},\mathbb{C})$ and
$H^t(M^{r,\bullet}(\mathcal{U}_\mathcal{K}),\delta)$ is zero if $r\neq q$ or $t \neq q-p.$ 
The pairing between $H_{-p+2q}(\mathbb{C}^n\setminus Z_{\mathcal{K}},\mathbb{C})$ and  $H^{-p+2q}(M^{\bullet}(\mathcal{U}_\mathcal{K}),D)$ is nondegenerate, therefore, 
the pairing between $H_{-p,2q}(\mathbb{C}^n\setminus Z_{\mathcal{K}},\mathbb{C})$ and $H^{q-p}(M^{q,\bullet}(\mathcal{U}_\mathcal{K}),\delta)$ is nondegenerate, and these spaces are dual. 
Hence, we have a natural isomorphism
$$H^{q-p}(M^{q,\bullet}(\mathcal{U}_\mathcal{K}),\delta) \stackrel{\phi\circ \varepsilon }{\simeq} H^{-p,2q}(R_\mathcal{K}\otimes\mathbb{C});$$
using this isomorphism and Lemma \ref{lemma.filtr} we obtain the statement of the theorem.
 \hfill$\square$


\end{document}